\documentclass[12pt]{article}

\usepackage[leqno]{amsmath}
\usepackage{amsthm,amscd}
\usepackage{amsfonts} 
\usepackage{amssymb} 

\newtheorem{theorem}{Theorem}[section]
\newtheorem{prop}[theorem]{Proposition}
\newtheorem{co}[theorem]{Corollary}
\newtheorem{lm}[theorem]{Lemma}
\newtheorem{rem}[theorem]{Remark}

\def\Proof{\noindent{\sl Proof.}\qquad}
\def\QED{\hfill\hbox{\vrule width 4pt height 6pt depth 1.5pt}\par\bigskip}

\def\diag{\mathop{\rm diag}}

\def\rank{\mathop{\rm rank}}

\newcommand{\trans}{^T}

\def\kmms{\kern-\mathsurround}

\newcommand{\calD}{{\cal D}}

\newcommand{\calA}{{\cal A}}

\newcommand{\bbR}{{\mathbb R}}

\newcommand{\bbZ}{{\mathbb Z}}
\newcommand{\bbQ}{{\mathbb Q}}

\newcommand{\period}{{\rho}}

\begin{document}

\title{Periodicity of hyperplane arrangements 
with integral coefficients modulo positive integers}
\author{Hidehiko Kamiya
\footnote
{
{\it Faculty of Economics, Okayama University}
}\\
Akimichi Takemura 
\footnote
{\it Graduate School of Information Science and Technology,
University of Tokyo}
\\
Hiroaki Terao 
\footnote
{
{\rm This work was
supported by the
MEXT
and the JSPS.}\,\,
{\tt hterao00@za3.so-net.ne.jp}\,\,\,
{\it Department of Mathematics, Hokkaido University}
}
}
\date{\today}
\maketitle
\begin{abstract}
  We study central hyperplane arrangements with integral coefficients
  modulo positive integers $q$.  We prove that the cardinality of the
  complement of the hyperplanes is a quasi-polynomial in two ways, first
  via the theory of elementary divisors and then via the theory of the 
  Ehrhart quasi-polynomials.  This result is useful for determining the
characteristic polynomial of the corresponding real arrangement.
With the former approach, we also prove that intersection lattices
modulo $q$ are periodic except for a finite number of $q$'s.  

\end{abstract}

\section{Introduction}
\label{sec:intro}

%

When a linear form in 
$x_{1}, \dots, x_{m}$ 
with integral coefficients
is given, we may naturally consider its ``$q$-reduction''
for any positive integer $q$.
The $q$-reduction 
is the image by the modulo $q$ projection
$
\bbZ[x_{1}, \dots, x_{m}]
\longrightarrow
\bbZ_{q} [x_{1}, \dots, x_{m}],
$ 
where $\bbZ_{q} = \bbZ /q \bbZ$. 
In this paper, we call the kernel of the
resulting linear form
a ``hyperplane'' in $V := \bbZ_{q}^{m}$.
Suppose that a finite set of nonzero linear forms with
integral coefficients is given. 
Then it not only defines a central
hyperplane arrangement $\calA$ in $\bbR^{m}$,  
but also gives
a ``hyperplane arrangement'' $\calA_{q}$ 
in $V$ through the $q$-reduction
for each $q\in \bbZ_{> 0} $.
%
%
%
A basic fact we prove in this paper is that
the cardinality of the complement 
$M(\calA_{q})$ of the arrangement $\calA_{q}$ in $V$,
as a
function of $q$, is a quasi-polynomial in $q$.  
(In other words, there exist a positive integer
$\rho$ (a period) and polynomials 
$P_{j} (t) \,\,(1\le j \le \rho)$ such that
$|M(\calA_{q} )| = P_{r} (q)
\,\,
(1\le r\le \rho,\,\,
q \in r + \rho\bbZ_{\ge 0})$ 
for all $q\in \bbZ_{> 0}$.) 
We provide two proofs of
this fact.  The first proof uses the theory of elementary divisors.  
The second proof is based on the
theory of the Ehrhart quasi-polynomials applied to each chamber of the
arrangement.  

In our setting, 
the approach via elementary divisors is more powerful 
than the one via the Ehrhart theory.  
The former gives more
information on the coefficients of the quasi-polynomials, 
and it also enables us to prove 
that the intersection lattices modulo $q$ are themselves periodic
except for a finite number of $q$'s.  
Despite the advantage of the
approach via elementary divisors for our setting, we also consider the
connection to the Ehrhart theory an important aspect of our discussion,
because many results in the Ehrhart theory can be applied to further
develop the arguments 
in this paper.

Especially when $q$ is a prime, the arrangement $\calA_{q}$
lies in the vector space $V = \bbZ_{q}^{m}$.
In this case, it is well known (e.g., \cite{crr}, \cite[(4.10)]{tertohoku},
\cite[Thm.3.2]{kott})
that 
$|M(\calA_q)|$ is equal to 
$\chi(\calA_{q}, q)$ 
and that
$\chi(\calA_{q}, t)$ 
coincides with
$\chi(\calA, t)$ 
for a sufficiently large prime $q$,
where $\chi( - , t)$ stands for 
 the characteristic polynomial
(e.g., \cite[Def.2.52]{ort},
\cite[Chap.3, Ex.56]{sta})
of an arrangement.
These facts provide the ``finite field method'' to study
the real arrangement $\calA$. 
The method was 
initiated and systematically applied by
%
%
Athanasiadis \cite{ath96, ath99, ath04}. 
It has been used to solve problems related to hyperplane arrangements
by
%
Bj\"orner and Ekedahl \cite{bje} and Blass and Sagan \cite{bls}
among others. 
It was also used in \cite{kott} to find
the characteristic polynomials of the mid-hyperplane
arrangements up to a certain
dimension.
Athanasiadis \cite{ath06} studies a problem similar to but different from the problem 
in the present paper. 
He proves that the coefficients of the characteristic polynomial of a certain 
deformation of a central arrangement are quasi-polynomials. 
A series of works by Athanasiadis on the finite field method is worth special mention 
as the driving force of the research on this method. 

For the theory of hyperplane arrangement, the reader is referred to 
\cite{ort}. 
For the Ehrhart theory for 
counting lattice points in rational polytopes, see the book by Beck and 
Robins \cite{ber}. 
Beck and Zaslavsky \cite{bez} study 
the extension of the Ehrhart theory to counting lattice points 
in ``inside-out polytopes''.

The organization of the paper is as follows.  In the rest of this
section, we set up our notation.  
In Section \ref{sec:quasi}, we prove that 
the cardinality of the complement
$M(\calA_{q})$ 
 is a quasi-polynomial in $q$, via the theory of
elementary divisors (Section \ref{subsec:via-elementary-divisors})
and via the theory of the Ehrhart quasi-polynomials (Section
\ref{subsec:ehrhart}).  
Based on this result, we consider a way of calculating the characteristic 
polynomial 
$\chi(\calA, t)$ of the corresponding real arrangement 
$\calA$ (Section \ref{subsection:characteristic-polynomial}).
In Section \ref{sec:intersection-lattice}, 
we prove that the intersection lattices modulo $q$ are
periodic except for a finite number of $q$'s.

In our forthcoming paper
\cite{ktt}, 
we apply the results in the present paper to 
the arrangements arising from root systems
and the mid-hyperplane arrangements. 



\subsection{Setup and notation}
\label{subsec:setup}
Let $m, n \in \bbZ_{>0}$ be positive integers. In this paper, $m$ denotes the dimension
and $n$ is the number of hyperplanes in an arrangement.  
Suppose we are given 
an $m \times n$ integer matrix 
\begin{equation*}
C=(c_1,\ldots,c_n) \in {\rm Mat}_{m \times n}(\bbZ)
\end{equation*}
consisting of 
column vectors 
$c_j=(c_{1j},\ldots,c_{mj})\trans \in \bbZ^m, \ 1\le j\le n$. 
Here, $\trans$ denotes the transpose and 
${\rm Mat}_{m \times n}(\bbZ)$ stands for the set of
$m\times n$ matrices with integer elements.
We assume that integral vectors 
$c_j$
are nonzero: 
\begin{equation}
\label{eq:cj-nonzero}
c_j\ne (0,\ldots,0)\trans, \quad 1\le j\le n. 
\end{equation}

Consider a real central hyperplane arrangement 
\[
\calA=\calA_C:=\{H_j: 1\le j\le n\}
\]
with 
\[
H_j=H_{c_j}:=\{x=(x_1,\ldots,x_m) \in \bbR^m: x c_j=0\}. 
\]
%
As an example, let us take $m=2, \ n=3$ and 
\begin{equation}
\label{eq:c1c2c3}
C=
\begin{pmatrix}
1 & 1 & -2 \\ 
-1 & 1 & 1 
\end{pmatrix},
\end{equation} 
i.e., $c_1=(1,-1)\trans, \ c_2=(1,1)\trans, \ c_3=(-2,1)\trans$. 
Then the corresponding 
hyperplane arrangement in 
$\bbR^2=\{ (x,y): x,y\in \bbR \}$ 
is ${\cal A}=\{ H_1, H_2, H_3\}$ with 
\[
H_1: x-y=0, \quad 
H_2: x+y=0, \quad
H_3: -2x + y=0.
\]

Since the coefficient vectors $c_j=(c_{1j},\ldots,c_{mj})\trans\in \bbZ^m, \ 1\le j\le n$, 
defining 
$H_j
$ are integral, 
we can consider the 
reductions 
of $c_j$ 
modulo positive integers 
$q\in \bbZ_{>0}$. 
Fix $q\in \bbZ_{>0}$ and 
let 
\[
[c_j]_q=([c_{1j}]_q,\ldots,[c_{mj}]_q)\trans\in \bbZ_q^m
\] 
be the {\it $q$-reduction} of 
$c_j $, 
i.e., 
$[c_{ij}]_q=c_{ij}+q\bbZ \in \bbZ_q, 
\ 1\le i \le m, \ 1\le j \le n$. 
In $V=\bbZ_q^m$, 
let us consider 
\[ 
H_{j, q}=H_{c_j, q}:=\{ x=(x_1,\ldots,x_m) \in V: x[c_j]_q=[0]_q \}, 
\]
and define 
\[
\calA_q=\calA_{C, q}:=\{ H_{j, q}: 1\le j \le n \}. 
\] 
We emphasize that $\calA_q=\calA_{C, q}$ is determined by 
$C$ and $q$, but not by $\calA=\calA_C$ and $q$. 
For a non-prime $q$, it may not be appropriate to call 
$H_{j, q}$ a hyperplane, but by abusing the terminology we call 
$H_{j, q}$ a hyperplane, and 
$\calA_q$ 
an arrangement of hyperplanes. 
In our previous example \eqref{eq:c1c2c3}, 
$\calA_q=\{ H_{1,q}, H_{2,q}, H_{3,q} \}$ 
with 
\begin{align}
H_{1, q} &=  \{ ([0]_q,[0]_q), ([1]_q,[1]_q), \dots, ([q-1]_q,[q-1]_q) \}, \nonumber \\
H_{2, q} &=  \{ ([0]_q,[0]_q), ([1]_q,[q-1]_q), \ldots, ([q-1]_q,[1]_q) \}, \label{eq:H2q-example} \\
H_{3, q} &=  \{ ([0]_q,[0]_q), ([1]_q,[2]_q), ([2]_q,[4]_q), \ldots, ([q-1]_q,[q-2]_q) \}. \nonumber 
\end{align}

%

In the finite field method and its generalization in the present paper, 
we are interested in the cardinality of the 
complement of $\calA_q$.   We denote the complement by 
\begin{equation*}
\label{eq:Mq}
M({\calA}_q) := V \setminus \bigcup_{1\le j \le n}H_{j, q}
\end{equation*}
and its cardinality by $|M({\calA}_q)|$.  
We will prove
that
$|M({\calA}_q)|$ is a quasi-polynomial in $q$ 
of degree $m$ and with the
leading coefficient identically equal to 1. 
That is, there exist 
{\it a period} $\period \in \bbZ_{>0}$ and 
$\alpha_{h, s}\in \bbQ, \ 0\le h\le m-1, \ s\in \bbZ_\period$, such that  
\begin{equation}
\label{eq:quasi-polynomial}
|M({\calA}_q)| 
= q^m + \alpha_{m-1,[q]_\period}q^{m-1} + \dots + \alpha_{1,[q]_\period}q 
+ \alpha_{0,[q]_\period},
\quad
q \in \bbZ_{>0}; 
\end{equation}
in fact, $\alpha_{h,s}, \ 0\le h\le m-1, \ s\in \bbZ_\period$, 
are integral: $\alpha_{h,s}\in \bbZ$. 
In this paper, we will call \eqref{eq:quasi-polynomial} 
the {\it characteristic quasi-polynomial} 
of $\calA_q$, 
because, as we will see in Section \ref{subsection:characteristic-polynomial}, 
the value \eqref{eq:quasi-polynomial} coincides with  
$\chi(\calA, q)$
if $q$ and $\rho$ are coprime, where $\chi(\calA, t)$ denotes 
the characteristic polynomial 
(e.g., \cite[Def.2.52]{ort},
\cite[Chap.3, Ex.56]{sta})
of the real arrangement 
$\calA$. 
The minimum period is simply called
{\it the period} of $|M({\calA}_q)|$.  
Often it is not trivial to find the period of 
$|M({\calA}_q)|$, although it is relatively easy to evaluate some
multiple of the period, which we simply call a period.
%

This is because of the following. 
The sum $\chi_1(q)+\chi_2(q)$ of two quasi-polynomials $\chi_1(q), \chi_2(q)$ 
is a quasi-polynomial having as a period the least common multiple
of the periods of $\chi_1(q)$ and $\chi_2(q)$.  
However, due to possible cancellations of terms, 
the period of $\chi_1(q)+\chi_2(q)$ may be smaller than this least common multiple. 
See McAllister and Woods \cite{mcw}. 

For a subset $J=\{j_1,\dots,j_k\} \subseteq \{ 1,\ldots,n\}$, write
\begin{equation}
\label{eq:HJq}
H_{J, q}:=\bigcap_{j\in J}H_{j, q}
= H_{j_1, q} \cap \cdots \cap H_{j_k, q}.
\end{equation}
When $J$ is nonempty, 
$H_{J, q}$ in \eqref{eq:HJq} is determined by the $q$-reduction of 
the $m\times k$ submatrix
\[
C_{J}:=(c_{j_1},\ldots,c_{j_k}) 
\in {\rm Mat}_{m \times k}(\bbZ) 
\]
of $C$; 
when $J$ is empty, we understand that $H_{\emptyset, q}=V$. 

The Smith normal form of an integer matrix 
$G\in {\rm Mat}_{m\times k}(\bbZ), \ k\in \bbZ_{>0}$, is 
\begin{eqnarray}
\label{eq:snfC}
\qquad 
SGT=
\begin{pmatrix}
E & O \\ 
O & O 
\end{pmatrix}
\in {\rm Mat}_{m \times k}(\bbZ), 
&& 
E=\diag(e_1,\ldots,e_\ell), \ \ 
\ell=\rank G, 
\\ 
&& 
\quad 
e_1,\ldots,e_{\ell}\in \bbZ_{>0}, \ \ e_1|e_2|\cdots|e_{\ell}, \nonumber 
\end{eqnarray}
where $S\in {\rm Mat}_{m\times m}(\bbZ)$ and $T\in {\rm Mat}_{k\times k}(\bbZ)$ 
are unimodular matrices. 
The positive integers $e_1,\ldots,e_{\ell}$ are 
the {\it elementary divisors} of $G$. 
For 
simplicity, we often use the following notation
\begin{equation*}
\label{eq:notation-diag}
\diag(\{e_1, \dots,e_\ell\}; m,k)
=
\begin{pmatrix}
E & O \\ 
O & O 
\end{pmatrix}
\in {\rm Mat}_{m \times k}(\bbZ) .
\end{equation*}

\section{Characteristic quasi-polynomial}
\label{sec:quasi}

\subsection{Via elementary divisors}
\label{subsec:via-elementary-divisors}

In this subsection, 
we prove that 
$|M({\calA}_q) |=|V \setminus \bigcup_{1\le j \le n}H_{j, q}|$ 
is a quasi-polynomial in  $q\in \bbZ_{>0}$ using the theory of 
elementary divisors. 

Let $I_Y( \, \cdot \, ), \ Y\subseteq V$, stand for the characteristic function 
(indicator function) of $Y: I_Y(x)=1, \ x\in Y$ and 
$I_Y(x)=0, \ x \in V\setminus Y$. 
Then for every $x \in V$, 
\[
\prod_{j=1}^n\left(1- I_{H_{j, q}}(x)\right) 
= \sum_{J \subseteq \{ 1,\ldots, n\}}(-1)^{|J|}I_{H_{J, q}}(x) 
= I_V(x)+\sum_{\emptyset \ne J \subseteq \{ 1,\ldots,n\}}
(-1)^{|J|}I_{H_{J, q}}(x), 
\]
which may be viewed as the inclusion-exclusion principle. 
%
Therefore, from the relation 
$x \in 
M({\calA}_q)
\Leftrightarrow 
1=\prod_{j=1}^n(1-I_{H_{j, q}}(x))$,  
we have
\begin{eqnarray}
|M({\calA}_q)|
&=& \sum_{x\in V}\prod_{j=1}^n\left( 1-I_{H_{j, q}}(x)\right) \label{eq:|M|=SP}
= q^m+\sum_{\emptyset \ne J \subseteq \{ 1,\ldots,n\}}
(-1)^{|J|}\left| H_{J, q}\right|. 
\end{eqnarray}  
Hence it suffices to verify that 
for each nonempty subset $J=\{j_1,\ldots,j_k\}$
of $\{1,\ldots,n\}$, 
%
the cardinality $|H_{J, q} |$
is a quasi-polynomial in $q\in \bbZ_{>0}$. 
Actually, we can show that $|H_{J, q} |$ is a quasi-monomial 
with an integral coefficient. 

Fix 
$J=\{j_1,\ldots,j_k\}\ne \emptyset$ and 
consider 
$C_{J}=(c_{j_1},\ldots,c_{j_k}) 
\in {\rm Mat}_{m \times k}(\bbZ)$. 
For each $
q \in \bbZ_{>0}$, 
let us define $f_{J,q}:V=\bbZ_q^m \to \bbZ_q^k$ by 
\begin{equation}
\label{eq:x-mapsto}
x \mapsto x[C_J]_q, 
\end{equation}
where 
$[C_J]_q=([c_{j_1}]_q,\ldots,[c_{j_k}]_q) 
\in {\rm Mat}_{m \times k}(\bbZ_q)$ 
is the $q$-reduction of $C_J$.  
Then 
$|H_{J, q}|= |\ker f_{J,q}|$, so 
the problem reduces to proving that $|\ker f_{J,q}|$ is 
a quasi-monomial in $q$. 
This fact can be shown by using the following general lemma. 

\begin{lm}
\label{lm:homomorphism}
Let $m$ and $k$ be positive integers. 
Let $f: \bbZ^m\to \bbZ^k$ be a $\bbZ$-homomorphism. 
Then the cardinality of the kernel of the induced morphism 
$f_q: \bbZ_q^m\to \bbZ_q^k$ is a quasi-monomial of $q\in \bbZ_{>0}$. 
Furthermore, suppose $f$ is represented by 
a matrix $G \in {\rm Mat}_{m\times k}(\bbZ)$. 
Then this quasi-monomial $|\ker f_q|, \ q\in \bbZ_{>0}$, can be expressed as 
\begin{equation}
\label{eq:q^(m-l)C}
|\ker f_{q}|=(d_1(q)\cdots d_{\ell}(q))q^{m-\ell}, 
\end{equation}
where $\ell=\rank G$ and $d_j(q):={\rm gcd}\{e_j, q \}, \ 1\le j \le \ell$. 
Here, $e_1,\ldots,e_{\ell}\in \bbZ_{>0}, \ e_1|e_2|\cdots |e_{\ell}$, are the elementary 
divisors of $G$. 
In that case, the quasi-monomial $|\ker f_q|, \ q\in \bbZ_{>0}$, 
has the 
minimum
period $e_{\ell}$, where we consider $e_0$ to be one. 
\end{lm}

\Proof 
If $f$ is the zero $\bbZ$-homomorphism, then $|\ker f_q|=|\bbZ_q^m|=q^m$ and  
the theorem is trivially true. 
So we may assume that $f$ is not the zero $\bbZ$-homomorphism. 
Since $|\ker f_q|=q^m/|{\rm im}f_q|$, we will study $|{\rm im}f_q|$.  

Suppose $f$ is represented by an $m\times k$ integer 
matrix $G 
\in {\rm Mat}_{m\times k}(\bbZ)$. 
Then, for $q\in \bbZ_{>0}$, the induced morphism $f_q: \bbZ_q^m\to \bbZ_q^k$ is 
given by $x \mapsto x[G]_q$. 

Consider the Smith normal form of $G$ in \eqref{eq:snfC}.  
%
Since unimodularity is preserved under $q$-reductions, 
%
we may assume that $G$ is of the form  
\begin{equation*}
G= \diag(\{e_1,\dots,e_\ell\}; m,k)
\end{equation*}
from the outset. 
Then we have  
\[
f_{q}(x)=([e_1]_q x_1,\ldots,[e_{\ell}]_q x_{\ell},
[0]_q,\ldots,[0]_q) \in \bbZ_q^k
\] 
for $x=(x_1,\ldots,x_m)\in \bbZ_q^m$.  
Therefore, 
${\rm im} f_{q}=[e_1]_q \bbZ_q\times \cdots \times [e_{\ell}]_q \bbZ_q$ 
and hence 
\[
|{\rm im} f_{q}|=\frac{q}{d_1(q)}\times \cdots \times \frac{q}{d_{\ell}(q)}
=\frac{q^{\ell}}{d_1(q)\cdots d_{\ell}(q)},
\] 
where 
$d_j(q)={\rm gcd}\{e_j, q \}, \ 1\le j \le \ell$. 
Consequently, we obtain \eqref{eq:q^(m-l)C}. 

Now, for any $j=1,\ldots,\ell$, we have 
$
d_j(q + e_{\ell})
={\rm gcd}\{ e_j, q+e_{\ell}\}={\rm gcd}\{ e_j, q\}
=d_j(q)$.
Therefore, \eqref{eq:q^(m-l)C} is a quasi-monomial in $q$ 
of degree $m-\ell <m$ and with a period $e_{\ell}$.   
In fact, we can show that $e_{\ell}$ is the minimum period as follows. 

Let $e'$ be the minimum period. 
Note $e'|e_{\ell}$. 
%
We have
$d_{j}(e_{\ell})=e_{j}\ge d_{j}(e')=d_{j}(e'+e_{\ell})>0$ for all $j=1,
\ldots,\ell.$
Since $e'$ is a period,
$d_1(e_{\ell})\cdots d_{\ell}(e_{\ell})
=
d_1(e'+e_{\ell})
\cdots d_{\ell}(e'+e_{\ell})$.
Therefore $e_{\ell} = d_{\ell} (e_{\ell} )
=
d_{\ell}(e' + e_{\ell} )
=
e'$.
\QED


Now, $f_{J,q}:V=\bbZ_q^m \to \bbZ_q^k$ in \eqref{eq:x-mapsto} 
is induced from
the $\bbZ$-homomorphism $f_J: \bbZ^m \to \bbZ^k$ 
represented by $C_J$. 
Thus, Lemma \ref{lm:homomorphism} 
implies that 
\begin{equation}
\label{eq:|H|}
|H_{J,q}|=|\ker f_{J,q}|=
(d_{J,1}(q)\cdots d_{J,\ell(J)}(q))q^{m-\ell(J)}
\end{equation} 
is a quasi-monomial with the period $e_{J,\ell(J)}$, 
where 
$\ell(J):=\rank C_J$ and  
$d_{J,j}(q):={\rm gcd}\{e_{J,j}, q \}$, $1\le j \le \ell(J)$. 
Here, $e_{J,1},\ldots,e_{J,\ell(J)}\in \bbZ_{>0}, \ 
e_{J,1}|e_{J,2}|\cdots |e_{J,\ell(J)}$, denote the elementary 
divisors of $C_J$. 
Note that $\ell(J)>0$ for all $J, \ |J|\ge 1$, because of the assumption 
\eqref{eq:cj-nonzero}. 

\begin{rem}
Assume that $q$ is prime. 
Then each $d_{J,j}(q)={\rm gcd}\{ e_{J,j},q\}, \ 1\le j\le \ell(J),$ is 
$1$ or $q$, and $d_{J,j}(q)=q$ if and only if $[e_{J,j}]_q=0$. 
It follows from \eqref{eq:|H|} that $X:=H_{J,q}$ for any nonempty $J$ 
satisfies $|X|=q^{m-\ell'}=q^{{\rm dim}X}$, where 
$\ell'=| \{ j: 1\le j\le \ell(J), \ [e_{J,j}]_q\ne 0 \} |$.  
Note that $|X|=q^{{\rm dim}X}$ for $X=H_{J,q}$ is true also when 
$J$ is empty: $|\bbZ_q^m|=q^m$.  
\end{rem}

From the discussions so far, 
we reach the following conclusions.  
First, 
$|M({\cal A}_q)|$, 
$q\in \bbZ_{>0}$,  
is a monic quasi-polynomial in $q$ of degree $m$. 
Second, 
a period of this quasi-polynomial can be obtained in the following way.
For each $m\times k \ (1\le k \le n)$ submatrices $C_J$ of 
$C=(c_1,\ldots,c_n)\in {\rm Mat}_{m \times n}(\bbZ)$, 
find its largest 
elementary divisor $e(J):=e_{J,\ell(J)}$. Let
\begin{equation*}
\label{eq:period-0}
 \period_0:={\rm lcm}\{ e(J) : J \subseteq \{1,\ldots,n\}, J \neq \emptyset \}.
\end{equation*}
Then 
$\period_0$ is a period of 
$|M({\cal A}_q)|$. 

For computing $\period_0$ when $m<n$, 
we can restrict the size of $J$ as $|J|\le m$: 
\begin{equation}
\label{eq:period--0}
 \period_0 ={\rm lcm}\{ e(J) : J \subseteq \{1,\ldots,n\}, \ 
1\le |J| \le \min\{m, n\} \}.
\end{equation}
We can prove \eqref{eq:period--0} in the following way. 
First, we note the next 
%
%
lemma.

\begin{lm} 
\label{lm:restriction}
Let $f_1, f_2: \bbZ^n\to \bbZ^m$ be two $\bbZ$-homomorphism with 
$\rank ({\rm im}f_1)=\rank ({\rm im}f_2)$ and 
${\rm im} f_2\subseteq {\rm im} f_1$. 
Then the largest elementary divisor of $f_1$ divides the largest elementary 
divisor of $f_2$. 
\end{lm}

\Proof
Define $I_{i} = {\rm Ann}({\rm coker}f_{i})
:=\{ p\in \bbZ: p({\rm coker}f_{i} )=0\},
$ and
the ideal $I_{i} $ is generated
by the largest elementary divisor of $f_i
\,\,
(i = 1, 2).$
Since there is a natural projection ${\rm coker}f_2\to{\rm coker}f_1$, 
we have $I_{2} \subseteq I_{1}$. This shows the lemma.  
\QED

Now, suppose $m<n$, and take an arbitrary $J\subseteq \{ 1,\ldots, n\}$ 
with $m<|J|\le n$. 
Let $\ell=\ell(J)=\rank C_J\  (\le\! m)$.  
Then we can take a subset $\tilde{J}\subset J, \ |\tilde{J}|=\ell$, such that 
$\rank C_{\tilde{J}}=\ell$. 
For this $\tilde{J}$, we have  
${\rm im}g_{\tilde{J}}\subseteq {\rm im}g_J$, where  
$g_J, g_{\tilde{J}}: \bbZ^n\to \bbZ^{m}$ are the $\bbZ$-homomorphisms 
defined by $C_J$ and $C_{\tilde{J}}$, respectively: 
$g_J(x)=\sum_{j\in J}x_jc_j, \ g_{\tilde{J}}(x)=\sum_{j\in \tilde{J}}x_jc_j, 
\ x=(x_1,\ldots,x_n)\trans \in \bbZ^n$. 
Then Lemma \ref{lm:restriction} implies that $e(J)|e(\tilde{J})$. 
From this observation, we obtain \eqref{eq:period--0}. 
When $n$ is considerably larger than $m$, the restriction $|J|\le m$ 
is computationally very useful.

Let us find a period $\period_0$ for our example 
\eqref{eq:c1c2c3}. 
Take $J=\{1,2\}$. 
Then we have 
\[
C_J=
\begin{pmatrix}
1 & 1 \\ 
-1 & 1 
\end{pmatrix}
\] 
with the Smith normal form 
$\diag(1, 2)$. 
Hence $e(J)=e_{J,2}
=2$. 
In a similar manner, we can find $e(J)$ for the other 
$J$'s with $1\le |J| \le 2$,  
and obtain $\period_0={\rm lcm}\{ 1,1,1,2,1,3
\}=6$.

Furthermore, 
for $ |J| \ge 1, \ 1\le j\le \ell(J)$,  
\begin{equation}
\label{eq:d(q)===}
d_{J,j}(q)={\rm gcd}\{ e_{J,j}, q \}={\rm gcd}\{ e_{J,j}, \period_0, q \}
={\rm gcd}\{ e_{J,j}, {\rm gcd}\{ \period_0, q\} \}. 
\end{equation}
This implies that the coefficient $d_{J,1}(q)\cdots d_{J,\ell(J)}(q)$ 
of each monomial 
$|H_{J, q}|, \ 
|J|\ge 1$, in \eqref{eq:|H|} 
depends on $q$ only through 
$\gcd\{ \period_0, q\}$. 
Therefore, the constituents of 
the quasi-polynomial $|M(\calA_q)|$ in \eqref{eq:|M|=SP} 
coincide for all $q$ with the same $\gcd\{ \period_0, q\}$.

We summarize the results obtained so far as follows:

\begin{theorem}
\label{thm:quasi-polynomial-1} 
The function 
$|M({\calA}_q)|$ is a monic quasi-polynomial in $q\in \bbZ_{>0}$ of degree
$m$ with a period $\period_0$ given in \eqref{eq:period--0}.  
Furthermore, 
in (\ref{eq:quasi-polynomial})  with $\period=\period_0$, 
the coefficients 
$\alpha_{h,[q]_{\period_0}}\in\bbZ, \ 0\le h \le m-1$, depend on $[q]_{\period_0}$ 
only through $\gcd\{ \period_0, q\}$. 
\end{theorem}


Let us find the characteristic quasi-polynomial for our example \eqref{eq:c1c2c3}. 
Since $\period_0=6
$, we know by Theorem \ref{thm:quasi-polynomial-1} that 
each of the sets $\{ 1,5 \}, \{ 2,4 \}, \{ 3,9\},\{ 6,12\}$ of values of $q$ 
determines a constituent of the characteristic 
quasi-polynomial $|M(\calA_q)|$. 
For $q=1$, we have $V=H_{1,1}=H_{2,1}=H_{3,1}=\{ ([0]_1, [0]_1) \}$ and thus  
$|M(\calA_1)|=0$. 
For $q=5$, we 
can count $|H_{1,5}\cup H_{2,5}\cup H_{3,5}|=13$ 
and get  
$|M(\calA_5)|=5^2-13=12$. 
By interpolation, we obtain the constituent 
$q^2-3q+2$ for 
$\gcd\{ 6, q\}=1$. 
In this way, 
we can get the following characteristic quasi-polynomial: 
\begin{equation}
\label{eq:quasi-example}
| M(\calA_q)|= 
\begin{cases}
q^2-3q+2 & \text{ when } {\rm gcd}\{ 6, q \}=1, \\ 
q^2-3q+3 & \text{ when } {\rm gcd}\{ 6, q \}=2, \\ 
q^2-3q+4 & \text{ when } {\rm gcd}\{ 6, q \}=3, \\ 
q^2-3q+5 & \text{ when } {\rm gcd}\{ 6, q \}=6. 
\end{cases} 
\end{equation}
From this characteristic 
quasi-polynomial, we can see that the minimum period is $6=\period_0$.  



\subsection{Via the Ehrhart theory}
\label{subsec:ehrhart}

We want to show via the Ehrhart theory that 
$|M({\cal A}_q)|=|V \setminus \bigcup_{1\le j \le n}H_{j, q}|$ 
is a quasi-polynomial in $q\in \bbZ_{>0}$.  The Ehrhart theory is indeed
useful for establishing that $|M({\cal A}_q)|$  is a
quasi-polynomial, and gives a geometric insight into its period. 
However, it does not seem to give information on the constituents of the
quasi-polynomial. 

For $j=1,\ldots,n$,  let 
\[
S_j:=\bbZ \cap \{ x c_j \mid x\in [0,1)^m \}.
\]
For example, for $c_j=(1,\dots,1)\trans \in \bbZ^m$
\[
\min_{x\in [0,1)^m \atop x_1+\cdots +x_m\in \bbZ} (x_1 + \dots + x_m)=0, \quad
\max_{x\in [0,1)^m \atop x_1+\cdots +x_m\in \bbZ} (x_1 + \dots + x_m)=m-1,
\]
and 
$S_j=\{ 0,1,\ldots,m-1 \}$ for this $c_j$.
%
Now define the additional  ``translated'' hyperplanes 
\[
H_{j}^{s_j}(q)=
H_{c_j}^{s_j}(q):=
\{ x=(x_1,\ldots,x_m) \in \bbR^m: x c_j=s_jq \}\subset \bbR^m, 
\quad s_j \in S_j,  
\]
for $j=1,\ldots,n,$ and consider the real hyperplane arrangement  
\[
\calA^{{\rm deform}}(q)=
\calA^{{\rm deform}}_C(q) := 
\{ H_{j}^{s_j}(q): 
s_j \in S_j, \ 1\le j \le n \}. 
\]
%
For any positive integer $q$, 
we can express 
$|M({\cal A}_q)|=|V \setminus \bigcup_{1\le j \le n}H_{j, q}|$ 
as 
\begin{eqnarray}
~~|M({\cal A}_q)|
&=&\left| \bbZ^m \cap [0,q)^m\setminus \bigcup 
\calA^{{\rm deform}}(q)\right| \label{eq:sharp-sharp} 
=\left| \bbZ^m \cap 
\left( q \times
\left([0,1)^m \setminus \bigcup 
\calA^{{\rm deform}}
\right)\right)\right|, 
\end{eqnarray}
where $\bigcup\calA^{{\rm deform}}(q):=\bigcup_{H\in \calA^{{\rm deform}}(q)}H$ 
and 
$\calA^{{\rm deform}}:=\calA^{{\rm deform}}(1)$.

Now, let us consider 
$[0,1)^m \setminus \bigcup \calA^{{\rm deform}}$ in \eqref{eq:sharp-sharp}.  
We see that $[0,1)^m$ is cut by  the hyperplanes 
$H_j^{s_j}(1), \ s_j \in S_j, \ 1\le j \le n$, 
into 
\[
P^O(s_1,\ldots,s_n):=
\{ x \in [0,1)^m: s_j < x c_j < s_j+1, \ 1\le j \le n \}, 
\]
$(s_1,\ldots,s_n) \in S_1^* \times \cdots \times S_n^*=:S^*$, 
where 
$S_j^*:=S_j\cup \{ \min S_j -1 \}, \ 1\le j\le n$. 
Therefore 
\begin{equation}
\label{eq:=sqcup}
[0,1)^m \setminus \bigcup \calA^{{\rm deform}}=
\bigsqcup_{(s_1,\ldots,s_n) \in S^*}P^O(s_1,\ldots,s_n) 
\end{equation}
is a disjoint union.
From \eqref{eq:sharp-sharp} and \eqref{eq:=sqcup}, we obtain 
\begin{eqnarray*}
|M({\calA}_q)|
&=&
\sum_{(s_1,\ldots,s_n) \in S^*}
\left| \bbZ^m \cap qP^O(s_1,\ldots,s_n) \right|
= \sum_{(s_1,\ldots,s_n) \in S^*} 
i(P^O(s_1,\ldots,s_n),q), 
\end{eqnarray*}
where $i(P^O(s_1,\ldots,s_n),q):=
|\bbZ^m \cap qP^O(s_1,\ldots,s_n)|.$

It should be noted that $P^O(s_1,\ldots,s_n), \ (s_1,\ldots,s_n) \in S^*$, 
are not necessarily open in $\bbR^m$.  
However, by applying the Ehrhart theory to some faces of each nonempty
$P^O(s_1,\ldots,s_n)$, we can show that
$i(P^O(s_1,\ldots,s_n),q)$
is a quasi-polynomial of $q\in \bbZ_{>0}$ with degree 
$
d
=\dim(P^O(s_1,\ldots,s_n))
$ 
and 
the leading coefficient equal to the normalized volume of 
$P^O(s_1,\ldots,s_n)$.
When $d=m$, the normalized volume is the 
same as the usual volume in $\bbR^m$. 
Therefore, we can conclude that the sum 
$\sum_{(s_1,\ldots,s_n) \in S^*} i(P^O(s_1,\ldots,s_n),q)=|M({\calA}_q)|$
is a quasi-polynomial 
of $q\in \bbZ_{>0}$ with 
degree 
$m$ 
and the leading coefficient 
$
\sum
{\rm vol}_m(P^O(s_1,\ldots,s_n))=
{\rm vol}_m([0,1)^m)=1$, where
${\rm vol}_m( \, \cdot \, )$ denotes the usual volume in $\bbR^m$.

\smallskip
Let us move on to the investigation into 
periods of the characteristic quasi-polynomial 
$|M({\calA}_q)|, \ q\in \bbZ_{>0}$. 
{}From the above discussion, we see that 
a common multiple of periods of 
$i(P^O(s_1,\ldots,s_n), q), \ (s_1,\ldots,s_n) \in S^*$, 
is a period of 
$|M({\calA}_q)|$.
Let $\bar{P}(s_1,\ldots,s_n)$ denote the closure of 
$P^O(s_1,\ldots,s_n)$.  Define
the {\it denominator} $\calD(\calA^{{\rm deform}})$ of $\calA^{{\rm
    deform}}$ by
\begin{eqnarray*}
\calD(\calA^{{\rm deform}}) &:=& \min \{ q \in \bbZ_{>0}: 
\text{ all } 
q\bar{P}(s_1,\ldots,s_n), 
\ (s_1,\ldots,s_n) \in S^*, \\ 
&& \qquad \text{ are integral polytopes}
\}.
\end{eqnarray*} 
The Ehrhart theory now implies that
$\calD(\calA^{{\rm deform}})$  is a period of $|M(\calA_q)|$.

Put $\tilde{C}:=(C, I_m)$, where 
$I_m$ is the $m\times m$ identity matrix. 
For $J\subseteq \{ 1,\ldots,n+m\}$, let 
$\tilde{C}_{J}$ denote the $m \times |J|$ 
submatrix of $\tilde{C}$ consisting of the columns  
corresponding to the elements of $J$. 
Then, 
in view of Cramer's formula, we see that  
$\calD(\calA^{{\rm deform}})$ divides 
\begin{eqnarray*}
&& \period_{\rm E}
:= {\rm lcm}\{ \det(\tilde{C}_{J}): 
J\subset \{1,\ldots,n+m\} \\ 
&& \qquad \qquad \qquad 
\text{ \ such that \ } 
|J|=m \text{ \ and \ } 
\det(\tilde{C}_{J}) \ne 0 \}. \nonumber 
\end{eqnarray*}
Hence, the minimum period divides 
$\period_{\rm E}$. 
%
For $J$ with $|J\cap \{ n+1, \dots,n+m\}|=h<m$,  
the determinant $\det(\tilde{C}_{J})$ equals 
an $(m-h)\times (m-h)$ minor 
of $C$ up to sign. 
Therefore, we can also write
\begin{equation}
\label{eq:period-E}
\period_{\rm E} 
= {\rm lcm}\{ \text{nonzero $j$-minors of $C$}, \ 
1\le j\le m \}.
\end{equation}
Now, recall the well-known fact that 
$\bar{e}(J):=e_{J,1}\times \cdots \times e_{J,\ell(J)}$ 
is equal to the greatest common 
divisor of all the (nonzero) $\ell(J)$-minors of $C_J, \ J\ne \emptyset$, 
and note the relation $e(J)=e_{J, \ell(J)}|\bar{e}(J), \ J\ne \emptyset$. 
Then we can easily see from \eqref{eq:period--0} and \eqref{eq:period-E} 
that $\period_0 | \period_{\rm E}$. 
Therefore, $\period_0$ gives a tighter 
bound for the period of the characteristic quasi-polynomial 
$|M({\calA}_q)|$ than $\period_{\rm E}$.

In our working example \eqref{eq:c1c2c3}, we have 
$\period_{{\rm E}}={\rm lcm}\{ 1,1,-2,-1,1,1,2,-1,3\}=6$ and thus 
$\period_0=
\period_{{\rm E}}$.  
In general, if we obtain the characteristic quasi-polynomial by interpolation using 
$\period_{{\rm E}}$ as a period and find that $\period_{{\rm E}}$ happens to be 
the minimum period, then we know $\period_0=
\period_{{\rm E}}$.  

\subsection{Characteristic polynomial of the real arrangement}
\label{subsection:characteristic-polynomial}

Let $\chi(\calA, t)$ be the characteristic polynomial
of the 
real hyperplane arrangement 
$\calA=\{H_j: 1\le j\le n\}$, where  
$H_j=\{x\in \bbR^m: x c_j=0\}, \ 1\le j\le n$. 

\begin{theorem}
Let $\period$ be a period of the quasi-polynomial
$|M({\calA}_q)|$ and $q$ be a positive integer relatively
prime to $\rho$. 
Then 
$|M({\calA}_q)| = \chi(\calA, q)$.
\end{theorem}

\Proof
Choose $c\in \bbZ_{\ge 0} $ and
$q'\in \bbZ$ such that 
$q = \rho c + q'$ and
$1\le q'\le \rho$. 
By Theorem \ref{thm:quasi-polynomial-1},
there exist integers 
$\alpha_{0}, \dots, \alpha_{m-1}$ such that
$|M({\calA}_k)|
=
k^m+\alpha_{m-1} k^{m-1}+\cdots+\alpha_0
$ for 
all $k \in q' +\period\bbZ_{\ge 0} $. 
Since $q'$ and $\period$ are 
relatively prime, then by 
Dirichlet's theorem on arithmetic progressions (e.g., \cite{ser}), 
$q'+\period \bbZ_{\ge 0}$
contains an infinite number of primes. 
On the other hand, 
it is well known 
(e.g., \cite{crr} \cite[(4.10)]{tertohoku} \cite[Theorem 3.2]{kott})
that, 
when $k$ is a sufficiently large prime, 
$
|M({\calA}_k)|
$
coincides with $\chi(\calA, k)$. 
Remember that the characteristic polynomial 
$\chi(\calA, t)$ 
is a 
monic polynomial of degree 
$\dim(\bbR^m)=m$.  
This implies that
$
\chi(\calA, t) = 
t^m+\alpha_{m-1} t^{m-1}+\cdots+\alpha_0
$
and thus
$
\chi(\calA, q) = 
q^m+\alpha_{m-1} q^{m-1}+\cdots+\alpha_0
=
|M({\calA}_q)|.
$
\QED

The argument above implies that we can obtain the characteristic polynomial 
$\chi(\calA, t)$ by 
counting 
$|M(\calA_{q_i})|=
|\bbZ_{q_i}^m\setminus \bigcup_{1\le j\le n}H_{j, q_i}|
$
for an arbitrary set of $m$ distinct values 
$q_1,\ldots,q_{m}$
with
${\rm gcd}\{\period, \, q_i \}=1
\,\,\,
(1\le i\le m)$. 
Note that $q_1,\ldots,q_m$ need not be prime.  

When $q'$ and $\period$ are not relatively prime, 
$q'+\period\bbZ_{\ge 0}$ contains at most 
one prime (and this prime is not necessarily
``sufficiently large''), 
so the above argument does not hold. 
For $q'+\period\bbZ_{\ge 0}$ with such $q'$'s, 
we obtain different polynomials 
than $\chi(\calA, t)$.  
%

In our example \eqref{eq:c1c2c3}, the constituent 
of the characteristic 
quasi-polynomial \eqref{eq:quasi-example} 
for $1+6\bbZ_{\ge 0}$ and $5+6\bbZ_{\ge 0}$, i.e.,  
$\gcd\{ 6, q\}=1$,  
is the characteristic polynomial of 
$\calA$. 
Thus $\chi(\calA, t)=t^2-3t+2=(t-1)(t-2)$. 


\section{Periodicity of the Intersection lattice}
\label{sec:intersection-lattice}
In this section, we show that 
the intersection lattice $L_q=L(\calA_q)$ 
(e.g., \cite[2.1]{ort}, \cite[Chap.3, Ex.56]{sta} )
is periodic for large enough $q$. 
Let us begin with our working example to illustrate the
periodicity of the intersection lattice $L_{q}$.

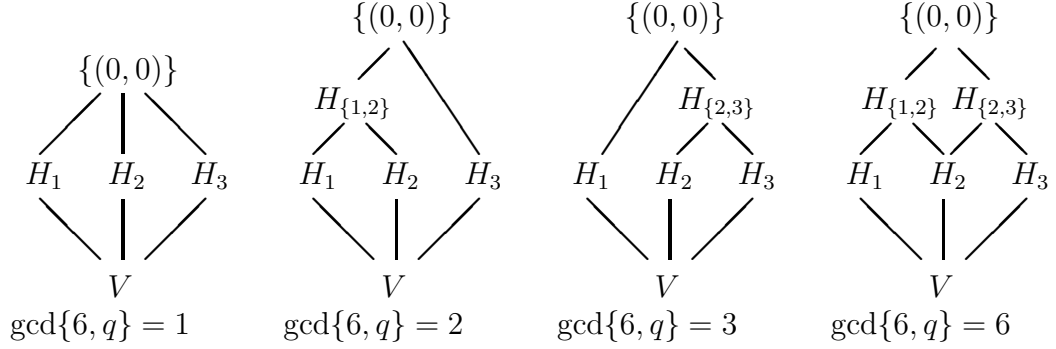
\begin{figure}
\begin{center}
\setlength{\unitlength}{1mm}
\begin{picture}(35,50)(0,-8)
\thicklines
\put(10,-5){${\rm gcd}\{ 6, q \}=1$}
\put(23,0){$V$}
\put(22,5){\line(-1,1){8}}
\put(25,5){\line(0,1){8}}
\put(28,5){\line(1,1){8}}
\put(12,15){$H_1$}
\put(23,15){$H_2$}
\put(34,15){$H_3$}
\put(14,19){\line(1,1){8}}
\put(25,19){\line(0,1){8}}
\put(36,19){\line(-1,1){8}}
\put(19,29){$\{(0,0)\}$}
\end{picture}
\begin{picture}(35,50)(0,-8)
\thicklines
\put(10,-5){${\rm gcd}\{ 6, q \}=2$}
\put(23,0){$V$}
\put(22,5){\line(-1,1){8}}
\put(25,5){\line(0,1){8}}
\put(28,5){\line(1,1){8}}
\put(12,15){$H_1$}
\put(23,15){$H_2$}
\put(34,15){$H_3$}
\put(14,19){\line(1,1){4}}
\put(25,19){\line(-1,1){4}}
\put(14,25){$H_{\{1,2\}}$}
\put(36,19){\line(-2,3){10}}
\put(20,29){\line(1,1){4}}
\put(19,36){$\{(0,0)\}$}
\end{picture}
\begin{picture}(35,50)(0,-8)
\thicklines
\put(10,-5){${\rm gcd}\{ 6, q \}=3$}
\put(23,0){$V$}
\put(22,5){\line(-1,1){8}}
\put(25,5){\line(0,1){8}}
\put(28,5){\line(1,1){8}}
\put(12,15){$H_1$}
\put(23,15){$H_2$}
\put(34,15){$H_3$}
\put(26,19){\line(1,1){4}}
\put(36,19){\line(-1,1){4}}
\put(26,25){$H_{\{2,3\}}$}
\put(16,19){\line(2,3){10}}
\put(31,29){\line(-1,1){4}}
\put(19,36){$\{(0,0)\}$}
\end{picture}
\begin{picture}(35,50)(0,-8)
\thicklines
\put(10,-5){${\rm gcd}\{ 6, q \}=6$}
\put(23,0){$V$}
\put(22,5){\line(-1,1){8}}
\put(25,5){\line(0,1){8}}
\put(28,5){\line(1,1){8}}
\put(12,15){$H_1$}
\put(23,15){$H_2$}
\put(34,15){$H_3$}
\put(26,19){\line(1,1){4}}
\put(36,19){\line(-1,1){4}}
\put(26,25){$H_{\{2,3\}}$}
\put(20,29){\line(1,1){4}}
\put(14,19){\line(1,1){4}}
\put(25,19){\line(-1,1){4}}
\put(14,25){$H_{\{1,2\}}$}
\put(31,29){\line(-1,1){4}}
\put(19,36){$\{(0,0)\}$}
\end{picture}
\caption{Hasse diagrams of intersection lattices for $q\ge 4$}
\label{fig:1}
\end{center}
\end{figure}

In our example 
\eqref{eq:c1c2c3},  
the ``hyperplanes'' $H_{1,q}, H_{2,q}, H_{3,q}$ were given in \eqref{eq:H2q-example}.  
For $q=1$, $V=\{([0]_1, [0]_1)\}=H_{1,1}=H_{2, 1}=H_{3, 1}$.  
For $q=2$, $H_{1, 2}=H_{2, 2}$ and 
for $q=3$, $H_{2, 3}=H_{3, 3}$.  These are the exceptions.
{}From  $q=4$ on, we have the periodicity of the
intersection lattice---the lattice of 
$H_{J,q}=\cap_{j\in J}H_{j,q}, \ J\subseteq \{ 1,2,3\}$, 
by reverse inclusion.
First, it is easily seen that 
for $q \ge 4$, $H_{j, q}$, $j=1,2,3$, are distinct, proper subsets of $V$.
Furthermore, for $q \ge 4$, 
\begin{eqnarray*}
H_{\{1,2\}, q}
&=&\begin{cases}
  \{([0]_q,[0]_q)\}, & q:\text{odd}, \\
  \{([0]_q,[0]_q), ([\frac{q}{2}]_q, [\frac{q}{2}]_q)\}, & q:\text{even},
\end{cases} \\ 
H_{\{2,3\}, q}
&=&\begin{cases}
  \{([0]_q,[0]_q)\}, & 3 \not| \ q , \\
  \{([0]_q,[0]_q), ([\frac{q}{3}]_q,[\frac{2q}{3}]_q), ([\frac{2q}{3}]_q,[\frac{q}{3}]_q)\}, & 3\;|\;q ,
\end{cases}
\end{eqnarray*}
and
$H_{\{ 1,3\},q}=H_{\{ 1,2,3\}, q}
=\{([0]_q,[0]_q)\}$.  We see that the intersection lattice for this example is
periodic and has  the period 6 for $q\ge 4$.  The Hasse diagrams for the four types of the 
intersection lattices are illustrated in Figure \ref{fig:1}.  In Figure
\ref{fig:1},  the subscript ${}_q$ is omitted for simplicity.

Let 
$J=\{ j_1,\ldots, j_k \}, \ 
1\le j_1<\cdots <j_k \le n, 
\ 1\le k \le n,$ 
be a nonempty subset of $\{ 1,\ldots,n \}$. 
We write the Smith normal form of $C_J\in {\rm Mat}_{m \times k}(\bbZ)$ as 
\begin{eqnarray}
\label{eq:SmithCJ}
&& S_JC_JT_J=
\diag(\{ e_{J,1},\ldots,e_{J,\ell(J)} \}; m, k)
=:\tilde{E}_J, \\ 
&& \qquad 
\ell(J)=\rank C_J, \quad 
e_{J,1},\ldots,e_{J,\ell(J)}\in \bbZ_{>0}, \quad  
e_{J,1}|e_{J,2}|\cdots|e_{J,\ell(J)}. \nonumber 
\end{eqnarray}
%
As in Section \ref{subsec:via-elementary-divisors}, 
we write $e(J)=e_{J,\ell(J)}$, 
the largest 
elementary divisor of $C_J$. 

Let $\period_0$ be 
the least common multiple of all $e(J)$'s with $1\le |J|\leq \min\{ m, n\}$ 
as in \eqref{eq:period--0}. 
Furthermore, 
define 
\[
q_0:=\max_{\emptyset \ne J\subseteq \{ 1,\ldots, n \}}
\min_{S_J} \, \max \{ | u |: u \text{ is an entry of } S_JC \text{ or } C \}, 
\] 
where the minimization is over all 
possible choices of $S_J$ in \eqref{eq:SmithCJ} for each fixed $J$. 



We are now in a position to state the main theorem of this section. 

\begin{theorem}
\label{th:period-L_q}
Let $J$ be an arbitrary nonempty subset of $\{ 1,\ldots, n\}$. 
Suppose $q, q' \in \bbZ_{>0}$ satisfy $q, q'>q_0$ and 
${\rm gcd}\{ \period_0, q \}={\rm gcd}\{ \period_0, q' \}$. 
Then, for any $j \in \{ 1,\ldots, n\}$, we have that $H_{j, q}\supseteq H_{J, q}$ 
if and only if $H_{j, q'}\supseteq H_{J, q'}$. 
\end{theorem}

When $j\in J$, the theorem is trivially true. 

In proving Theorem \ref{th:period-L_q}, we need the following 
proposition. 
Regard $V=\bbZ_q^m$ as a $\bbZ_q$-module. 
Let $V^*$ be the $\bbZ_q$-module consisting of the linear forms on $V$. 
For any $A\subseteq V$, we denote by $I(A)$ the set of linear forms 
vanishing on $A$: 
\[
I(A)=\{ \alpha \in V^*: \alpha(x)=0 \text{ \ for all \ } x \in A \}. 
\] 
Also, for any $B\subseteq V^*$, let $V(B)$ stand for the set of points 
at which each linear form in $B$ vanishes: 
\[
V(B)=\{ x \in V: \alpha(x)=0 \text{ \ for all \ } \alpha \in B \}. 
\]  
Evidently, $I(A)$ and $V(B)$ are submodules of $V^*$ and $V$, respectively.   

\begin{prop}
\label{prop:I(V(B))}
For any $B\subseteq V^*$, we have $I(V(B))=\langle B\rangle$, where 
$\langle B \rangle$ signifies the submodule of $V^*$ spanned by $B$. 
\end{prop}

\Proof
It suffices to show $I(V(B))=B$ for any submodule $B$ of $V^*$. 
It is trivially true that $I(V(B))\supseteq B$, so we will prove $I(V(B))\subseteq B$.   

Let $C_q(B)\in {\rm Mat}_{m\times k}(\bbZ_q), \ k=|B|,$ be the 
coefficient matrix of $B$. 
We can find an integral matrix $C(B)\in {\rm Mat}_{m\times k}(\bbZ)$ 
whose $q$-reduction is $C_q(B)$, i.e., 
$[C(B)]_q=C_q(B)$. 
Now, let $e_1|e_2|\cdots | e_{\ell}, \ \ell=\rank C(B)$, 
be the elementary divisors of $C(B)$. 
In $\bbZ_q$, we then have $C_q(B)$ is equivalent to 
\begin{equation}
\label{eq:diag[]0}
\diag(\{[e_1],\ldots,[e_{\ell'}]\}; m,k)
\in {\rm Mat}_{m\times k}(\bbZ_q), 
\end{equation} 
where $[e_1], \ldots, [e_{\ell'}] \in \bbZ_q\setminus\{ 0 \}, \ \ell'\le \ell$. 
Here, we are writing $[ \, \cdot \, ]$ for $[ \, \cdot \, ]_q$ for simplicity. 
We can choose $C(B)$ in such a way that $\ell'=\ell$, and we decide to do so. 
From \eqref{eq:diag[]0} we see that we can assume $B$ is spanned by 
$[e_1]y_1, \ldots, [e_{\ell}]y_{\ell}$ after a suitable coordinate change, 
where $\{ y_1,\ldots,y_{\ell}, y_{\ell +1}, \ldots, y_m \}$ is a basis of $V^*$.  
It follows that $V(B)$ is spanned by 
\begin{align*}
p_1&:=(\,[q/d_1(q)]\,, 0,\ldots,0), \ldots,\;
p_{\ell}:=
 ( \underbrace{0,\ldots,0}_{\ell-1},\,[q/d_{\ell}(q)]\,,
 0,\ldots,0),  \\  
& p_{\ell+1}:=( \underbrace{0,\ldots,0}_{\ell}, 1, 0,\ldots,0 )
,\ldots, \;
p_m:=\left( 0,\ldots,0, 1 \right) 
\end{align*}
with $d_j(q)={\rm gcd}\{ e_j,q\}, \ 1\le j \le \ell$. 

Now, take an arbitrary $\alpha=[a_1]y_1+\cdots+[a_{m}]y_{m} 
\in I(V(B))=I(p_1,\ldots,p_m)$ 
with $[a_1],\ldots,[a_{m}]\in \bbZ_q$. 
Then we have $0=\alpha(p_1)=[qa_1/d_1(q)],$ so 
$a_1=r_1d_1(q)$ for some $r_1\in \bbZ$.  
This implies $[a_1]=[r_1][d_1(q)]=[r'_1][e_1]$ with 
$[r'_1]:=[r_1][e_1/d_1(q)]^{-1}\in \bbZ_q$, 
where $[e_1/d_1(q)]^{-1}$ exists because 
${\rm gcd}\{e_1/d_1(q), q \}=1$. 
Similarly, for each $j=2,\ldots,\ell$, we have $[a_j]=[r'_j][e_j]$ for some 
$[r'_j]\in \bbZ_q$. 
Moreover, for $j=\ell+1,\ldots,m$, we obtain $0=\alpha(p_j)=[a_j]$.  
Therefore, we have 
$\alpha = [r'_1][e_1]y_1+\cdots+ [r'_{\ell}][e_{\ell}]y_{\ell} \in B$, and 
the proof is complete.  
\QED
\begin{flushleft}
{\it Proof of Theorem \ref{th:period-L_q}.}
\end{flushleft}
Without loss of generality, we may assume $j=1$. 
Let $[S_J]_q, [C_J]_q, [T_J]_q
$ and $[\tilde{E}_J]_q$ be 
the $q$-reductions of 
$S_J, C_J, T_J
$ and $\tilde{E}_J$ in 
\eqref{eq:SmithCJ}, respectively. 

First, we know by Proposition \ref{prop:I(V(B))} 
that $H_{1,q}\supseteq H_{J, q}$ if and only if 
$[c_1]_q$ lies in the column space of 
$[C_J]_q$ in $\bbZ_q^m$. 
Since 
$S_J^{-1}$ and $T_J^{-1}$ exist in ${\rm Mat}_{m\times m}(\bbZ)$ and 
${\rm Mat}_{k\times k}(\bbZ)$, respectively,  
the latter condition is equivalent to 
$[c_1]_q$ being in the column space of $[C_J]_q[T_J]_q=
[S_J^{-1}]_q[\tilde{E}_J]_q$, 
which in turn is equivalent to  
$[S_J]_q[c_1]_q$ being in the column space of $[\tilde{E}_J]_q$ 
in $\bbZ_q^m$. 

Next, let us 
paraphrase the above condition in $\bbZ_q^m$ 
as a condition in $\bbZ^m$. 
The condition holds if and only if 
$S_Jc_1\in \bbZ^m$ is in the column space of 
$(\tilde{E}_J, qI_m)\in {\rm Mat}_{m\times (k+m)}(\bbZ)$ in $\bbZ^m$.
Noting that 
$e_{J,j}\bbZ+q\bbZ=d_{J,j}(q)\bbZ$ with  
$d_{J,j}(q)={\rm gcd}\{ e_{J,j}, q \}, \ 1\le j\le \ell(J)$, 
we see that the condition holds if and only if  
$S_Jc_1$ is in the column space of 
$\diag(d_{J,1}(q),\ldots,d_{J,\ell(J)}(q),q,\ldots,q)\allowbreak
\in {\rm Mat}_{m\times m}(\bbZ)$. 
Since the absolute value of each entry of $S_Jc_1\in \bbZ^m$ is 
less than $q$, the condition is equivalent to  
$S_Jc_1$ being in the column space of 
\begin{equation}
\label{eq:hline}
\diag(\{d_{J,1}(q), \dots, d_{J,\ell(J)}(q)\}; m, \ell(J)) \ 
\in {\rm Mat}_{m\times \ell(J)}(\bbZ). 
\end{equation}

Now, since the absolute value of each entry of $S_Jc_1$ is 
less than $q'$ as well, 
the preceding argument holds true also for $q'$. 
Moreover, we see from \eqref{eq:d(q)===} that  
$d_{J,j}(q)=d_{J,j}(q')$ for $j=1,\ldots,\ell(J)$. 
Thus \eqref{eq:hline} remains the same when $q$ is 
replaced by $q'$. 
Therefore, we obtain the desired result. 
\QED


Our assumption \eqref{eq:cj-nonzero} implies that 
$H_{j, q} \not\supseteq H_{\emptyset, q}=V, \ 1\le j\le n$, for all $q>q_0$. 
From this observation and Theorem \ref{th:period-L_q}, 
it follows immediately that $L_q=L(\calA_q)$ for $q>q_0$ is 
periodic in $q$ with a period $\period_0$. 

\begin{co}
\label{co:periodicity-of-L}
The intersection lattice $L_q=L(\calA_q)$ is 
periodic in $q>q_0$ with a period $\period_0$: 
\[
L_{q+s\period_0}\simeq L_q \text{ for all } q > q_0 \text{ and } s\in \bbZ_{\ge 0}. 
\]  
\end{co}

%
%

Finally, we make a remark on the coarseness of the intersection lattices
for different $q$'s. In Figure \ref{fig:1}
we see that the intersection lattice for the case $\gcd\{ 6, q\}=6$ is
the most detailed and that the coarseness is nested according to the
divisibility of $\gcd\{ 6, q\}$.
This observation can be generally stated as follows.

\begin{prop}
\label{prop:divisibility}
Let $I, J\subseteq \{1,\dots,n\}$ and suppose that
$H_{I,q}=H_{J,q}$  for some $q > q_0$.  Then 
$H_{I,q'}=H_{J,q'}$ for every $q'> q_0$ such that
$\gcd\{ \period_0, q' \}|\gcd\{ \period_0, q \}$.
\end{prop}

\Proof It suffices to show that for any $i\in I$, if $[c_i]_q$ lies in
the column space of $[C_J]_q$ in $\bbZ_q^m$, then
$[c_i]_{q'}$ lies in
the column space of $[C_{J}]_{q'}$ in $\bbZ_{q'}^m$. 
Without loss of generality, take $i=1$ and assume that 
$[c_1]_q$ lies in
the column space of $[C_{J}]_q$ in $\bbZ_q^m$. 
Then $S_J c_1$ is in the column space of \eqref{eq:hline}. 
Now, because $\gcd\{ \period_0, q'\}|\gcd\{ \period_0, q\}$ by assumption, 
we can see from \eqref{eq:d(q)===} that 
$d_{J,j}(q')|d_{J,j}(q), \ 1\le j\le \ell(J)$. 
This implies that 
$S_J c_1$ is in the column space of \eqref{eq:hline} with $q$ replaced by $q'$.  
Therefore, $[c_1]_{q'}$ lies in the column space of $[C_{J}]_{q'}$ in $\bbZ_{q'}^m$. 
\QED

\end{document}